\documentclass[12pt,centertags,oneside]{amsart}

\usepackage{amsmath,amstext,amsthm,amscd,typearea,hyperref,stmaryrd}
\usepackage{amssymb}
\usepackage{a4wide}
\usepackage[mathscr]{eucal}
\usepackage{mathrsfs}
\usepackage{typearea}
\usepackage{charter}
\usepackage{pdfsync}
\usepackage[a4paper,width=16.2cm,top=3cm,bottom=3cm]{geometry}

\numberwithin{equation}{section}

\newtheorem{theorem}{Theorem}[section]

\newtheorem{proposition}[theorem]{Proposition}
\newtheorem{corollary}[theorem]{Corollary}
\newtheorem{lemma}[theorem]{Lemma}

\newcommand{\cali}[1]{\mathscr{#1}}

\newcommand{\dist}{{\rm dist}}

\newcommand{\ddc}{{dd^c}}

\newcommand{\ddbar}{{\partial\overline\partial}}

\newcommand{\id}{{\rm id}}

\renewcommand{\Re}{{\rm Re}}

\newcommand{\Cc}{\cali{C}}

\newcommand{\Ec}{\cali{E}}

\renewcommand{\Mc}{\cali{M}}

\newcommand{\Pc}{\cali{P}}

\newcommand{\FS}{{\rm FS}}

\newcommand{\C}{\mathbb{C}}
\newcommand{\D}{\mathbb{D}}

\newcommand{\N}{\mathbb{N}}

\newcommand{\R}{\mathbb{R}}

\renewcommand{\S}{\mathbb{S}}
\renewcommand{\P}{\mathbb{P}}

\newcommand{\Ebf}{\mathbf{E}}

\newcommand{\Mf}{{\mathfrak{M}}}

\newcommand{\sign}{\mathop{\mathrm{sign}}\nolimits}
\newcommand{\Res}{\mathop{\mathrm{Res}}\nolimits}
\newcommand{\Sym}{\mathop{\mathrm{Sym}}\nolimits}
\newcommand{\Prob}{\mathop{\mathrm{Prob}}\nolimits}

\title[Large deviation theorem for random covariance matrices]{Large deviation theorem for random covariance matrices}

\author{Tien-Cuong Dinh}
\address{Department of Mathematics, National University 
of Singapore, 10 Lower Kent Ridge Road, Singapore 119076. 
{\tt Email: matdtc@nus.edu.sg ; http://www.math.nus.edu.sg/$\sim$matdtc} }

\author{Duc-Viet Vu}
\address{UPMC Univ Paris 06, UMR 7586, Institut de Math\'ematiques de Jussieu, 4 place Jussieu,
F-75005 Paris, France. {\tt Email: duc-viet.vu@imj-prg.fr}}

\date{July 19, 2017}

\begin{document}

\maketitle

\begin{abstract}
We establish a large deviation theorem for the empirical spectral distribution of random covariance matrices whose entries are independent random variables with mean 0, variance 1 and having controlled forth moments. Some new properties of Laguerre polynomials are also given. 
\end{abstract}

\medskip\medskip

\noindent
{\bf MSC 2010:} 60B20, 33C45. 

 \medskip

\noindent
{\bf Keywords:} random covariance matrices, Marchenko-Pastur law, Laguerre polynomials, large deviations.

\section{Introduction} \label{Intro}

Let $M=[t_{ij}]_{1 \le i \le p, 1 \le j \le n}$ be a random matrix  whose entries are independent random variables of mean $0$ and variance $1.$   The number $p$ may depend on $n.$
 The matrix  $W:= M^* M$ is called {\it a random covariance matrix} and is one of the most important random matrices used in statistical inference. Recall that $M^*$ is the conjugate transpose of $M$. 
 We are interested in the distribution of the eigenvalues of $n^{-1}W$ as $n$ and $p$ tend to infinity.
The matrices  $M^* M$ and $M M^*$ have essentially the same spectrum : their only difference
 is the multiplicities at 0 when $p\not=n$, see 
 Lemma \ref{lemma_singvalue} below. Therefore, without loss of generality, 
 we only consider the case where $p\le n.$ 
 
 Let $\lambda_1,\ldots,\lambda_n$ be the eigenvalues of $n^{-1}W$. Then {\it the empirical spectral distribution} of $n^{-1}W$ is the probability measure 
 $$\mu_{p,n}:={1\over n} \sum_{k=1}^n\delta_{\lambda_k}, $$
 where $\delta_{\lambda_k}$ denotes the Dirac mass at $\lambda_k$.
 Let $\phi \in (0,1]$,  $a:= (1 - \sqrt{\phi})^2$ and $b:= (1 + \sqrt{\phi})^2.$  Recall that the probability measure associated to {\it the Marchenko-Pastur law} is 
$$\mu_\phi:=(1- \phi) \delta_0+ \frac{1}{2\pi t} \sqrt{(t-a)(b-t)_{+}}\, dt,$$
where  $x_+:=\max(x,0)$.
This measure is supported by $\{0\}\cup [a,b]$ and contains an atom at 0 (except when $\phi=1$). Its restriction to $[a,b]$ has the density $\frac{1}{2\pi t} \sqrt{(t-a)(b-t)_{+}}$ with respect to the Lebesgue measure on $\R$.

If the entries $t_{ij}$ of $M$ are i.d.d.,  when $p,n$ tend to infinity and $p/n$ converges to a positive number $\phi_0$, then with probability 1,
the empirical spectral distribution of $n^{-1}W$ satisfies the well-known Marchenko-Pastur law.
That is, with probability 1, we have $\mu_{p,n}\to\mu_{\phi_0}$ as $p,n\to\infty$ and $p/n\to\phi_0$. For non i.i.d., the same property holds under a supplementary hypothesis on the tails of $t_{ij}$ which is automatically true in the i.i.d. case.
For the proof of the Marchenko-Pastur law, see Bai-Silverstein \cite{BS},  Marchenko-Pastur \cite{MP}, Wachter \cite{Wachter} and Yin \cite{Yin}. We also refer the reader to the following papers and the references therein for related results : Bai-Hu-Zhou \cite{BHZ},
 Bai-Yin \cite{BY}, Ben Arous-P\'ech\'e \cite{BP}, Bloemendal et al. \cite{BEK},
 Cacciapuoti-Maltsev-Schlein \cite{CMS},  G\"otze-Tikhomirov \cite{GT}, Grenander-Silverstein \cite{GS}, Jonsson \cite{Jonsson} and Tao-Vu \cite{TaoVu}.

In the present article, we assume that the following property  holds for all entries $t_{ij}$ 
\begin{align*} 
\Ebf(|t_{ij}|^4) \le \beta,
\end{align*}
where $\beta>0$ is a real number. The entries $t_{ij}$ are not supposed to be identical.
Here is our main result.

\begin{theorem} \label{th_main}
Let $M, W, t_{ij}, p,n,\mu_{p,n}$ and $\beta$ be as above.  Define $\phi:=p/n$. Then there are  universal constants $A_1>0$ and $A_2>0$  with the following property. For every $\delta>0$, 
there is a set  $\Ec_{p,n}(\delta)$ of $(p\times n)$-matrices satisfying the following estimate of probability
$$\Prob(M\in \Ec_{p,n}(\delta))\leq A_1n^{A_1} e^\beta e^{-A_2\delta n}$$
and such that if  $M\not\in \Ec_{p,n}(\delta)$ and $I\subset \R$ is an interval then
$$\dist(\mu_{p,n},\mu_\phi)\leq \delta \qquad \text{and} \qquad |\mu_{p,n}(I)-\mu_\phi(I)| \leq  {\sqrt{\delta}\over 1-\phi}\cdot$$
\end{theorem}
 
The distance $\dist(\cdot,\cdot)$ between probability measures will be introduced later. In the last theorem, it allows us to estimate the rate of convergence in the Machenko-Pastur law.
For instance, the Kantorovich-Wasserstein distance is bounded by a constant times $\dist(\cdot,\cdot)^{1/2}$, see Section \ref{section_Laguerre} for details. Note that the last inequality in the theorem is not useful when $\phi$ is close to 1. In this case, using some techniques from \cite{Dinh} we can obtain useful estimates. More precisely, a similar estimate holds when $I$ is outside a neighbourhood of $0$ (where the density of $\mu_\phi$ is big) and a weaker inequality holds for $I$ close to 0.  However, we will not consider this question here in order to keep the paper less technical.

Note also that $0$ is an eigenvalue of multiplicity at least $n-p$ of $n^{-1}W$, see Lemma \ref{lemma_singvalue} below. So $\mu_{p,n}$ has an atom at 0 of mass at least $1-\phi$. It follows that when $\phi\to 0$ we get $\mu_{p,n}\to\delta_0$ for any choice of $M$ and the rate of convergence depends on the rate of convergence of $\phi$ to 0. So  the only interesting case is when $\phi=p/n$ is bounded from below by a positive constant.

Since the above constants $A_1$ and $A_2$  do not depend on $p,n,\delta$ and $\beta$, one can apply the above result even when $\beta$ and $\delta$ depend on $n$. For example, if $\beta$ is a constant and $\phi$ converges to some number $\phi_0$, by taking suitable $\delta\gg n^{-1}\log n$, we get a rate for the almost sure convergence of $\mu_{p,n}$ to $\mu_{\phi_0}$, in terms of $\delta$ and in terms of the rate of convergence of $\phi$ to $\phi_0$. 

To prove the main result, we will use an abstract large deviation theorem for the distribution of the zeros of polynomials of degree $p$, see Theorem \ref{th_dinh} below. We will apply this theorem for the polynomial $z^{-(n-p)}\det(z-n^{-1}W)$ of degree $p$, whose zeros are essentially the eigenvalues of the matrix $n^{-1}W$. This approach requires an upper bound of the expectation of $|\det(z-n^{-1}W)|^2$, see Proposition \ref{prop_det_W} below, which will be obtained using a  long combinatoric computation in Section \ref{section_det}. We also need some properties of Laguerre's polynomials that will be presented in Section \ref{section_Laguerre}. Note that the computation in this section allows us to obtain new properties on the distribution of zeros of Laguerre's polynomials which are of independent interest, see Corollary \ref{cor_zeros_Laguerre} below.

\medskip
\noindent
{\bf Acknowledgement.} 
This work was supported by the Start-Up 
Grant R-146-000-204-133 from the National University of  Singapore.
It was partially written during the visit of the first author at the Freie Universit\"at Berlin
and of the second author at the National University of Singapore. They would like to thank these organisations, the Alexander von Humboldt foundation and H\'el\`ene Esnault for their hospitality and  support.

\section{Laguerre polynomials and Marchenko-Pastur law} \label{section_Laguerre}

In this section, we will give some estimates on Laguerre polynomials and also discuss the relation of these polynomials with the Marchenko-Pastur law. These properties will be used later in the proof of our main theorem. 
Let $\alpha$ be a real number. For a positive integer $p$  the function
$$L^{(\alpha)}_p(z):= \sum_{k=0}^p \binom{p+\alpha}{p-k} \frac{(-z)^{k}}{k!}$$
is called \emph{a generalized Laguerre polynomial}, where the generalized binomial coefficients are defined by
$$\binom{p+\alpha}{p}:={(p+\alpha)(p+\alpha-1)\ldots (\alpha+1)\over p!}\cdot$$
Here is the main estimate on Laguerre polynomials that we need in this paper.

\begin{proposition} \label{prop_Laguerre}
There is a universal constant $c>0$ such that we have for every $\alpha\geq 0$, $p\geq 1$ and $z\geq 0$
$$|L_p^{(\alpha)}((p+\alpha)z)|\leq c\min(p,1+p/\alpha) (p+\alpha)^{p/2}p^{-{p/ 2}} z^{-{\alpha/ 2}}e^{(p+\alpha)z/ 2} e^{-{\alpha/ 2}}.$$
\end{proposition}
\proof
We can assume $\alpha>0$ since the case where $\alpha=0$ can be obtained by continuity.
Define $n:=p+\alpha$, $\phi:=p/n$ and $\widetilde\phi:=\phi/(1-\phi)=p/\alpha$. We need to show  for some universal constant $c>0$ that
$$z^{{(n-p)/ 2}}e^{-nz/ 2}|L_p^{(n-p)}(nz)|\leq c\min(p,1+p/\alpha)  \phi^{-n\phi/2}  e^{-{(n-p)/ 2}}.$$
Observe that for $p\geq n/2$ we have
$$\min(p,1+\widetilde\phi)\geq {1\over 2} \min(n,1+\widetilde\phi)$$
and for $p\leq n/2$ we have $\alpha\geq n/2$, $\widetilde\phi\leq 1$ and hence
$$\min(p,1+\widetilde\phi)\geq 1 \geq {1\over 2} \min(n,1+\widetilde\phi).$$
So in both cases, for the desired inequality, we can (and we will do) replace  $\min(p,1+\widetilde\phi)$ by  $\min(n,1+\widetilde\phi)$.

Now, by using the integral formula for Laguerre polynomial \cite[p.105]{Szego}, we have
$$L_p^{(n-p)}(nz) ={1\over 2i\pi}\oint {e^{-{nzt\over 1-t}}\over (1-t)^{n-p+1}t^{p+1}} dt.$$
Here, the integral is taken on a simple counter-clockwise piecewise smooth contour about 0. Note that the same meaning will be used for the other integrals below.
Using a new variable $u=-t/(1-t)$ we have $t=-u/(1-u)$, $1-t=1/(1-u)$, $dt=-(1-u)^{-2}du$ and 
$$L_p^{(n-p)}(nz) =\pm{1\over 2i\pi}\oint \big[e^{2zu} (1-u)^2 u^{-2\phi}\big]^{n/2} du/u.$$
Thus,
$$z^{{(n-p)/2}}e^{-nz/2}L_p^{(n-p)}(nz) =\pm{1\over 2i\pi}\oint \Big[e^{-z}e^{2zu} z^{1-\phi}(1-u)^2 u^{-2\phi}\Big]^{n/2} du/u.$$

Write $u=re^{i\theta}$ and $\xi:=\cos\theta$.
We will consider a contour that will be specified later. It is contained in the half-plane $\{\Re(u)< 1/2\}$ and given by an equation $r=r(\xi)$.
We have
$$z^{{(n-p)/2}}e^{-nz/2}|L_p^{(n-p)}(nz)|\leq \int_{-\pi}^{\pi} 
\Big[e^{-(1-2r\xi)z} z^{1-\phi}(1+r^2- 2r\xi)r^{-2\phi}\Big]^{n/2} (1+|r'/r|)d\theta.$$
The expression in the brackets, seen as a function in $z\geq 0$, has its maximal value when $z=(1-\phi)/(1-2r\xi)$. By substituting this value of $z$ into the last integral, we deduce that this integral is bounded by
$$ e^{-(n-p)/2}\int_{-\pi}^{\pi} \Big[(1-\phi)^{1-\phi}(1-2r\xi)^{-(1-\phi)}(1+r^2- 2r\xi)r^{-2\phi}\Big]^{n/2} (1+|r'/r|)d\theta.$$
Denote by $g(r,\xi)$ the function in the last brackets with  $r\xi=\Re(u)<1/2$, $r>0$ and $|\xi|\leq 1$. 

\medskip\noindent
{\bf Case 1.} Consider first the case where $1+\widetilde\phi\leq n$ which implies that 
$$\min(n, 1+\widetilde\phi)\geq {1\over 2} (1+\widetilde\phi).$$ 
We will use  the contour defined by the following equivalent equations 
$$r^2+2\xi\widetilde\phi r - \widetilde\phi=0 \quad \Longleftrightarrow \quad r^2=\widetilde\phi(1-2r\xi) \quad \Longleftrightarrow \quad (1-2r\xi)=r^2/\widetilde\phi.$$
Since $\widetilde\phi>0$, for each $\xi\in[-1,1]$, these equations have a unique positive solution given by 
$$r=r(\xi):=-\widetilde\phi\xi +\sqrt{\widetilde\phi^2\xi^2+\widetilde\phi}={\widetilde\phi\over \widetilde\phi\xi+\sqrt{\widetilde\phi^2\xi^2+\widetilde\phi}} \geq{\widetilde\phi\over 2\widetilde\phi+\sqrt{\widetilde\phi}}\geq {1\over 3}\min(1,\sqrt{\widetilde\phi}).$$ 
It is clear from the above equivalent equations that the solution satisfies $1-2r\xi>0$ and hence $r\xi<1/2$. 
If we consider $g$ as a function in $r$ and $1-2r\xi$, it is not difficult to see that the contour is exactly the set where the differential of $g$ vanishes. 

A direct computation using $(1-2r\xi)=r^2/\widetilde\phi$ shows that $g(r,\xi)=\phi^{-\phi}$ on the considered contour. 
Therefore, to get the desired estimate, we only need to bound $1+|r'/r|$. From the above discussion, we have 
$$r\geq {1\over 3}\min(1,\sqrt{\widetilde\phi}) \quad \text{and}\quad |r'|=\Big| -\widetilde \phi+{\widetilde\phi^2\xi \over \sqrt{\widetilde \phi^2\xi^2+\widetilde\phi}}\Big|\leq 2\widetilde\phi.$$ 
So $1+|r'/r|$ is  bounded from above by a constant times $1+\widetilde\phi$ (we can easily see it by considering $\widetilde\phi\geq 1$ and $\widetilde\phi\leq 1$). The desired estimate follows.

\medskip\noindent
{\bf Case 2.} Consider now the case where $1+\widetilde\phi\geq n$. We have  $\min(n, 1+\widetilde\phi)=n$. We will use the contour  defined by 
$$r^2+2\xi n r - n=0 \quad \Longleftrightarrow \quad r^2=n(1-2r\xi) \quad \Longleftrightarrow \quad (1-2r\xi)=r^2/n.$$
As above, we obtain that $1+|r'/r|$ is bounded by a constant times $n$. Define 
$$\delta:=1-\phi=1/(1+\widetilde\phi)\leq 1/n.$$
A direct computation using $(1-2r\xi)=r^2/n$ gives us
$$g(r,\xi)^{n/2}=(1+1/n)^{n/2} \delta^{n\delta/2} n^{n\delta/2}\leq \sqrt{e}\leq \sqrt{e}\phi^{-n\phi/2}$$
since $\phi<1$. The proposition follows.
\endproof

It is well-known that the zeros of a Laguerre polynomial of suitable parameters are equidistributed with respect to the Marchenko-Pastur law when the degree of the polynomial tends to infinity, see e.g. Dette-Studden \cite{DetteStudden}. We will give at the end of this section the rate of this convergence, after recalling necessary notions and results.
We will first give some basic properties of probability measures on the complex plane $\C$, their logarithmic potentials, and 
 some notions of distance between these measures, see \cite[section 2]{Dinh} for details.

Let $\P^1=\C\cup\{\infty\}$ denote the Riemann sphere which is the natural compactification of $\C$ by adding a point $\infty$ at infinity. Let $z$ denote the standard complex coordinate in $\C$. Recall that the Fubini-Study form on $\P^1$ is defined by $\omega_\FS:=\ddc  \log (1+|z|^2)^{1/2}$ where the operator $\ddc:={i\over \pi} \ddbar$ can be identified to $1/(2\pi)$ times the Laplacian operator. The differential form $\omega_\FS$ extends to a smooth differential form on $\P^1$ and induces there a Hermitian metric that we will use here. 

For any positive measure $\mu$ with compact support in $\C$, its logarithmic potential $u$ is defined by 
$$u(z):=\int_\C \log|z-w| d\mu(w) \quad \text{for} \quad z\in\C.$$
This is the unique subharmonic function in $\C$ with values in $\R\cup\{-\infty\}$ such that if $m$ is the mass of $\mu$ then
$$\ddc u=\mu \quad \text{and} \quad \lim_{z\to \infty}u(z)-m\log|z|=0.$$
The first identity is understood in the sense of currents or distributions.

Let  $\Mc_c(\C)$ be the set of all probability measures with compact support in $\C$. 
For $\mu,\mu'$ in $\Mc_c(\C)$ and $u,u'$ their logarithmic potentials, consider the following notions of distance
$$\dist(\mu,\mu'):=\|u-u'\|_{L^1(\P^1)}:=\int_{\P^1} |u-u'| \omega_\FS$$
and for every $\gamma>0$
$$\dist_\gamma(\mu,\mu'):=\sup \Big\{ |\langle \mu-\mu',\phi\rangle|, \ \phi \text{ is a } \Cc^\gamma \text{ function on } \P^1 \text{ with } \|\phi\|_{\Cc^\gamma}\leq 1\Big\},$$
where the pairing $\langle \mu-\mu',\phi\rangle$ denotes the integral of $\phi$ with respect to the measure $\mu-\mu'$. 
Note that $\dist_1$ is equivalent to the well-known Kantorovich-Wasserstein distance. 
We have the following propositions, see \cite[section 2]{Dinh} and \cite{DinhSibony}.

\begin{proposition} \label{prop_dist_gamma}
{\rm (i)} For any $0<\gamma\leq \gamma'$, there is a constant $c_{\gamma,\gamma'}>0$ depending only on $\gamma$ and $\gamma'$ such that 
$$\dist_{\gamma'}\leq \dist_\gamma\leq c_{\gamma,\gamma'}[\dist_{\gamma'}]^{\gamma/\gamma'}.$$

{\rm (ii)}  For every $0<\gamma\leq 2$, there is a constant $c_\gamma>0$ depending only on $\gamma$ such that 
$$\dist_\gamma\leq c_\gamma \dist^{\gamma/2}.$$
\end{proposition}

\begin{proposition} \label{prop_dist_sup}
Let $L$ be a compact interval in the real line $\R$. Let $K\subset L$ be a compact interval  and let 
$\mu_0$ be a probability measure with support in $K$ whose
logarithmic potential $u_0$ is continuous. Consider  another probability measure $\mu$ with compact support in $\C$ and its logarithmic potential $u$.

{\rm (i)}  There is a constant $c>0$ depending only on $L$  such that
$$\dist(\mu,\mu_0) \leq c\sup_K (u-u_0).$$ 

{\rm (ii)} For every $0<\gamma\leq 2$, there is a constant $c_\gamma>0$ depending only on $L$ and $\gamma$ such that 
$$\dist_\gamma(\mu,\mu_0)\leq c_\gamma\sup_K (u-u_0)^{\gamma/2}.$$ 

{\rm (iii)} Assume that $\mu_0$ is absolutely continuous with respect to the Lebesgue measure on $\R$ and its density function is bounded by a constant $A$. Assume also that $\mu$ has support in $\R$. 
Then there is a constant $c>0$ depending only on $L$ such that for any interval $I\subset \R$
$$|\mu(I)-\mu_0(I)| \leq cA\big[\sup_K (u-u_0)\big]^{1/2}.$$
\end{proposition}
\proof
There are only two minor precisions in this statement when we compare it with the results in  \cite[section 2]{Dinh}. 
Firstly, the role of $A$ in the
last property (iii) is  can be easily seen from the proof in the above reference. Observe that this constant $A$ is bounded from below by the inverse of the length of $L$.
Secondly, in that reference, we used a conformal map $\phi:\P^1\setminus K\to \P^1\setminus \overline\D$ with $\phi(\infty)=\infty$, where $\D$ is the unit disc in $\C$. Also in this reference, we observed using the maximum principle that $\sup_K (u-u_0)=\sup_\C (u-u_0)$, see also the end of the proof of Corollary \ref{cor_Laguerre_MP}. Therefore, we can replace $K$ with $L$ and hence the map $\phi$ does not depend on $\mu_0$.
This is the reason why the constants involving in the proof depend only on $L$. 
\endproof

We will see that the distribution of the zeros of the above Laguerre polynomials is related to the Marchenko-Pastur law.
Let $\phi, a, b, \mu_\phi$ be as in the Introduction.  
Let  $\mu^+_\phi$ be the probability measure which is equal to $1/\phi$ times the non-atomic part of $\mu_\phi$, i.e.
$$\mu^+_\phi:= \frac{1}{2\pi \phi t} \sqrt{(t-a)(b-t)_{+}}\, dt.$$ 
Let $u_\phi$ and $u_\phi^+(z)$ denote respectively the logarithmic potentials of $\mu_\phi$ and $\mu_\phi^+$.
Since the logarithmic potential of $\delta_0$ is $\log|z|$, we have 
$$u_\phi(z)=(1-\phi) \log|z| + \phi u_\phi^+(z).$$

\begin{proposition} \label{prop_potential_MP}
We have for $z \in [a,b]$ 
$$u_\phi(z)=\frac{z-1-\phi}{2} +  \frac{1- \phi}{2} \log z + \frac{\phi\log \phi}{2}$$
and
$$u_\phi^+(z)= \frac{z-1-\phi}{2\phi} -  \frac{1- \phi}{2\phi} \log z + \frac{\log \phi}{2}\cdot$$
\end{proposition}
\proof 
Note that the first identity follows from the second one and the relation between $u_\phi$ and $u_\phi^+$ mentioned above. We prove now the second identity and
we only consider $\phi \in (0,1)$ since the case $\phi=1$ can be deduced by continuity. 

Observe that $b-a= 4 \sqrt{\phi}$ and $b+a= 2(1+\phi).$ Define a new variable $t'$ by
 $$t= \frac{1}{2} \frac{b-a}{2}t'+ \frac{a+b}{2}= \sqrt{\phi}t'+ (1+\phi).$$ 
It satisfies $t' \in [-2,2]$ and $dt= \sqrt{\phi} dt'$ when $t\in[a,b]$. Moreover,
$$(t-a)(b-t)=\Big[{1\over 4}(b-a)\Big]^2 (4-t'^2)=(4-t'^2)\phi.$$
Define also $z'$ by putting $z= \sqrt{\phi}z'+ (1+\phi)$, and  the function $u$ by $u(z'):= u^+_\phi(z).$ We have using the definition of the logarithmic potential $u_\phi^+(z)$
\begin{align*}
u(z') &= \int_a^b \log|z-t| {1\over 2\pi\phi t} \sqrt{(t-a)(b-t)} dt \\
&= \frac{1}{2\pi \phi} \int_{-2}^2 \log |z'- t'|  \frac{\sqrt{(4- t'^2)\phi}}{\sqrt{\phi}t'+ (1+\phi)}  \sqrt{\phi} dt' +  \frac{\log \sqrt{\phi}}{2\pi \phi} \int_{-2}^2 \frac{\sqrt{(4- t'^2)\phi}}{\sqrt{\phi}t'+ (1+\phi)}  \sqrt{\phi} dt' \\
&=\frac{1}{2\pi} \int_{-2}^2 \log |z'- t'|  \frac{\sqrt{4- t'^2}}{\sqrt{\phi}t'+ (1+\phi)} dt' + \frac{\log \sqrt{\phi}}{2\pi} \int_{-2}^2 \frac{\sqrt{4- t'^2}}{\sqrt{\phi}t'+ (1+\phi)} dt'. 
\end{align*}
We want to compute $u(z')$ for $z\in[a,b]$ or equivalently for $z'\in[-2,2]$. 

Let $u_1(z')$ and $u_2(z')$ denote respectively the first and second terms in the last sum. 
Note that the map $w \longmapsto w+ 1/w$ is 2 to 1 from the unit circle $\S^1$ to the interval $[-2,2]$. So we write 
$z'=w+1/w=w+\overline w$ with $w\in\S^1$, and
 use the new variable $s=e^{i\vartheta}\in\S^1$ such that $t'=s+1/s=2\cos\vartheta$ with $-\pi\leq \vartheta\leq 0$.
Since $\log|w|=\log|s|=0$, we have 
\begin{eqnarray*}
\log|z'-t'| & = & \log |1-we^{i\vartheta}|+\log |1-we^{-i\vartheta}| \\
& = & \Re \big[\log (1-we^{i\vartheta})+\log (1-we^{-i\vartheta})\big],
\end{eqnarray*}
where we use the principal branch for the complex logarithmic function. It follows that
\begin{eqnarray*} 
u_1(z') & = & \Re \Big[{1\over 2\pi}\int_{-\pi}^0 {\log (1-we^{i\vartheta}) 4\sin^2\vartheta d\vartheta\over \sqrt{\phi}(e^{i\vartheta}+e^{-i\vartheta})+ (1+\phi)}+{1\over 2\pi}\int_{-\pi}^0 {\log (1-we^{-i\vartheta}) 4\sin^2\vartheta d\vartheta\over \sqrt{\phi}(e^{i\vartheta}+e^{-i\vartheta})+ (1+\phi)} \Big] \\
& = & \Re\Big[{1\over 2\pi}\int_{-\pi}^\pi {\log (1-we^{i\vartheta}) 4\sin^2\vartheta d\vartheta\over \sqrt{\phi}(e^{i\vartheta}+e^{-i\vartheta})+(1+\phi)} \Big] \\
& = & \Re \Big[  - \frac{1}{2 \pi i} \int_{\S^1} \frac{\log(1 - w s)(s - 1/s)^2}{\sqrt{\phi}(s+ s^{-1})+ (1+\phi) } {ds\over s}\Big].
\end{eqnarray*}

Let $f(s)$ be the integrand in the last integral which is a meromorphic function in $s.$ We have 
$$f(s)=\frac{\log(1 - w s)(s-1/s)^2}{\sqrt{\phi}(s^2+ 1)+ s(1+\phi) }=\frac{\log(1 - w s)(s-1/s)^2}{\sqrt{\phi}(s+ \sqrt{\phi})(s+ 1/\sqrt{\phi})} \cdot$$
Since $\phi \in (0,1),$  the poles of $f$ in the unit disk are  $0$ and $-\sqrt{\phi}$. Both poles are simple because $\log(1-ws)=-ws+o(ws)$ vanishes when $s=0$.  Thus, by residue theorem (a nice logarithmic singularity in $\S^1$ doesn't cause any problem here), $u_1(z')$ is the real part of
$$- \frac{1}{2\pi i} \int_{\S^1} f(s) ds= -\Res(f,0)- \Res(f,-\sqrt{\phi})= \frac{w}{\sqrt{\phi}}- \frac{\log(1+ w \sqrt{\phi})}{ \sqrt{\phi}} (1/ \sqrt{\phi} - \sqrt{\phi}).$$
Here, we use the fact that if $a$ is a simple pole of a function $g(s)$ then $\Res(g,a)$ is equal to the value of $(s-a)g(s)$ when $s=a$.
So we obtain 
$$u_1(z')= \Re \Big[ \frac{w}{\sqrt{\phi}}- \frac{\log(1+ w \sqrt{\phi})}{ \sqrt{\phi}} (1/ \sqrt{\phi} - \sqrt{\phi}) \Big ].$$
 
In the same way, we get
$$u_2(z')= - \frac{\log\sqrt{\phi}}{4 \pi i} \int_{\S^1} \frac{(s-1/s)^2}{\sqrt{\phi}(s+ \sqrt{\phi})(s+1/\sqrt{\phi}) } ds.$$
Now, the point $-\sqrt{\phi}$ is a simple pole but 0 is a double pole. If $a$ is a double pole of a function $g(s)$ then $\Res(g,a)$ is equal to the derivative of $(s-a)^2g(s)$ at the point $a$. A direct computation gives us
$$u_2(z')=- \frac{\log \sqrt{\phi}}{2\sqrt{\phi}}\big [ - (1/ \sqrt{\phi} +\sqrt{\phi})  + (1/ \sqrt{\phi} - \sqrt{\phi}) \big ]=  \log \sqrt{\phi}.$$
Thus,
\begin{eqnarray*}
u(z') & = & \Re \Big [ \frac{w}{\sqrt{\phi}}- \frac{\log(1+ w \sqrt{\phi})}{ \sqrt{\phi}} (1/ \sqrt{\phi} - \sqrt{\phi}) \Big]+  \log \sqrt{\phi} \\
& = & \frac{w+ \bar{w}}{2 \sqrt{\phi}} - \frac{1- \phi}{2 \phi} \log[(1+ w\sqrt{\phi} )(1+ \bar{w}\sqrt{\phi})] +  \log \sqrt{\phi}\\
&=  & \frac{z'}{2\sqrt{\phi}} -  \frac{1- \phi}{2\phi} \log[1+\phi+ \sqrt{\phi} z']+ \frac{\log \phi}{2} \cdot
\end{eqnarray*}
Substituting $z'= (z-1-\phi)/ \sqrt{\phi}$ into the last equality gives
$$u_\phi^+(z)=\frac{z-1-\phi}{2 \phi} -  \frac{1- \phi}{2\phi} \log z + \frac{\log \phi}{2}\cdot$$
This ends the proof of the proposition.
\endproof

Let $a_1,\ldots,a_p$ denote the zeros of the rescaled Laguerre polynomial $L^{(\alpha)}_p((p+\alpha)z)$. It is known for $\alpha>0$ that all zeros of this Laguerre polynomial are real positive numbers, \cite[Th. 3.3.4]{Szego}. Denote by $\mu_p^{(\alpha)}$ the empirical measure for the zeros of  $L^{(\alpha)}_p((p+\alpha)z)$, i.e.
$$\mu_p^{(\alpha)}:={1\over p} (\delta_{a_1}+\cdots+\delta_{a_p}).$$
Observe also that  if $P(z)$ is a monic polynomial of degree $p$, then $p^{-1}\log|P(z)|$ is the logarithmic potential of the empirical measure of the zeros of $P(z)$. So the logarithmic potential of $\mu_p^{(\alpha)}$ is equal to 
$u(z):=p^{-1}\log |p!(p+\alpha)^{-p} L_p^{(\alpha)}((p+\alpha)z)|$ since the coefficient of $z^p$ in the polynomial $p!(p+\alpha)^{-p} L_p^{(\alpha)}((p+\alpha)z)$ is $\pm 1$.

\begin{corollary} \label{cor_Laguerre_MP}
Let $\alpha>0$ be a real number and $p\geq 1$ be an integer number. Define $\phi:=p/(p+\alpha)$.
Then there is a universal constant $c>0$ such that for any $z\in \C$
$$|p!(p+\alpha)^{-p} L_p^{(\alpha)}((p+\alpha)z)|\leq c\min(p,1+p/\alpha)p^{1/2}e^{pu^+_\phi(z)}.$$
\end{corollary}
\proof
We use the notation $n:=p+\alpha$ and $\phi=p/n$. 
We first consider the case where $z\in [a,b]$.
By Stirling's formula, we have $p!\lesssim p^{1/2} p^pe^{-p}$. So by Proposition \ref{prop_Laguerre}, the left hand side of the desired estimate is bounded by a constant times
$$p^{1/2}p^pe^{-p}n^{-p} \min(p,1+p/\alpha) n^{p/2}p^{-p/2}z^{-(n-p)/2}e^{nz/2}e^{-(n-p)/2}.$$
On the other hand, Proposition \ref{prop_potential_MP} implies that 
$$pu_\phi^+(z)={nz -n-p\over 2} - {n-p\over 2}\log z +{p\log(p/n)\over 2}\cdot$$
This proves the corollary for $z\in [a,b]$.

Now, let  $\delta>0$ be the real number  such that $e^{p\delta}=c\min(p,1+p/\alpha)p^{1/2}$ and consider $u(z)$ as above. 
The estimate obtained for $z\in [a,b]$ implies that $u-u_\phi^+\leq \delta$ on $[a,b]$.
Since $\mu_\phi^+$ is supported by $[a,b]$, its logarithmic potential $u_\phi^+$ is harmonic in $\C\setminus [a,b]$. It follows that $u-u_\phi^+$ is subharmonic in $\C\setminus [a,b]$. Moreover, using the properties of logarithmic potentials, we have $u(z)-u_\phi^+(z)\to 0$ as $z\to\infty$. So the function $u-u_\phi^+$ extends to a subharmonic function on $\P^1\setminus[a,b]$. We have seen that $u-u_\phi^+\leq \delta$ on $[a,b]$ which is the boundary of $\P^1\setminus[a,b]$. By maximum principle, the same inequality also holds on $\P^1\setminus [a,b]$. This implies the corollary for all $z\in\C$.
\endproof

\begin{corollary} \label{cor_zeros_Laguerre}
Let $\alpha>0$, $0<\gamma\leq 2$ be real numbers and $p\geq 2$ be an integer number. Define $\phi:=p/(p+\alpha)$.
Then there are a universal constant $c>0$  and a constant $c_\gamma>0$ depending only on $\gamma$ such that 
$$\dist(\mu_p^{(\alpha)},\mu_\phi^+)\leq cp^{-1}\log p \qquad \text{and} \qquad \dist_\gamma(\mu_p^{(\alpha)},\mu_\phi^+)\leq c_\gamma p^{-\gamma/2}(\log p)^{\gamma/2}.$$
Moreover, for any interval $I$ in $\R$, we have 
$$|\mu_p^{(\alpha)}(I)-\mu_\phi^+(I)| \leq c(1-\phi)^{-1} p^{-1/2} (\log p)^{1/2}.$$
\end{corollary}
\proof
We will use the notation given at the beginning of the section. Observe that $\mu_\phi^+$ is supported by $[a,b]$ which is a subset of $[0,4]$.  Moreover, its density with respect to the Lebesgue measure is bounded by a constant times
$$\sup_{t\geq a} {\sqrt{t-a}\over t}=\sup_{s\geq 0} {s\over s^2+a}={1\over 2\sqrt{a}}
={1\over 2(1-\sqrt{\phi})}\leq {1\over 1-\phi}\cdot$$
Using Corollary \ref{cor_Laguerre_MP}, we obtain the result as 
a direct consequence of Proposition \ref{prop_dist_sup} applied to the measures $\mu:=\mu_p^{(\alpha)}$ and 
$\mu_0:=\mu_\phi^+$. 
\endproof

Note that the last corollary improves  a well-known result by Dette-Studden which says that when $\alpha\to\infty$, $\phi\to\infty$ and $p/(p+\alpha)$ converges to some constant $\phi>0$, the measure $\mu_p^{(\alpha)}$ converges weakly to $\mu_\phi^+$, see  \cite{DIW,DetteStudden}. 
Note also that the last property in the corollary is not useful when $\phi$ tends very fast to 1. However, since the density of $\mu^+_\phi$ for $\phi=1$ is bounded outside a neighbourhood of 0, we can prove a similar property for $I$ outside a neighbourhood of $0$. Near the point 0, a weaker estimate can  also be obtained, see \cite{Dinh} for necessary techniques.

\section{Spectral distribution of random covariance matrices} \label{section_det}

In this section, we will give the proof our main theorem.
A key ingredient of the proof  is an upper bound for the expectation of $|\det(z-n^{-1}W)|^2$.
More precisely, we have the following result, where we use the notations introduced in the previous sections.

\begin{proposition} \label{prop_det_W}
There is a universal constant $c>0$ such that 
$$\Ebf \big(|\det(z-n^{-1}W)|^2\big)  \le c p^{13/2} n^{1/2}  e^{\beta} e^{2n u_\phi(z)}$$
for   $z \in [a,b]$ and  $n,p \in \N$ with $1\leq p\leq n$, and $\phi:=p/n$.
\end{proposition}

This estimate will be obtained as a consequence of a long combinatoric computation that will be presented below in a sequence of lemmas.  The following basic property has been used to reduce the study of our matrices to the case with $p\leq n$ by replacing $M$ with $M^*$.

\begin{lemma} \label{lemma_singvalue} 
We have $\det(z- M^* M)= z^{n-p} \det(z- M M^*).$
\end{lemma}

\proof By the singular value decomposition for $M,$ we have  $M= U \tilde{M} V^*,$ where $\tilde{M}$ is a real rectangular  diagonal  $(p\times n)$-matrix and $U,V$ are unitary matrices. It follows that $M^*M=V \tilde{M}^* \tilde{M} V^*$ and $M M^*=U \tilde{M} \tilde{M}^* U^*.$ Therefore, we get
$$\det  (z- M^* M)=\det(z- \tilde{M}^* \tilde{M})=z^{n-p}\det(z- \tilde{M} \tilde{M}^*)= z^{n-p} \det(z- M M^*).$$
The lemma follows.
\endproof

Write $W= [\xi_{jk}]_{1 \le j,k \le n}.$ Since $M=[t_{ij}]_{1\leq i\leq p,1\leq j\leq n}$ and $W=M^*M$, we have
$$\xi_{jk}:= \sum_{l=1}^p \bar{t}_{lj} t_{lk},$$
for $1 \le j,k \le n.$
Define $\llbracket1,n\rrbracket:=\{1,\ldots,n\}$. If $J$ is a set, denote by $\Sym(J)$ the symmetric group of all permutations of $J$, $\Mf(J)$ the set of all maps from $J$ to $\llbracket 1,p\rrbracket$, and $\Mf^\star(J)$ the set of injective maps from $J$ to $\llbracket 1,p\rrbracket$. 
For any set $J\subset\llbracket 1,n\rrbracket$, $\sigma \in \Sym(J)$ and $\tau\in\Mf(J)$ define 
\begin{align} \label{eq_T}
T(J,\sigma,\tau):=  \prod_{j \in J} \bar{t}_{\tau(j),j} t_{\tau(j), \sigma(j)}=   \prod_{j \in J} \bar{t}_{\tau(j),j} t_{\tau(\sigma^{-1}(j)), j}
\end{align}
and
\begin{align} \label{eq_xisigma}
\xi_{\sigma}:= \prod_{j \in J} \xi_{j \sigma(j)}= \sum_{\tau\in\Mf(J)}  T(J,\sigma,\tau)
\qquad \text{and} \qquad
\xi_{\sigma}^\star:= \sum_{\tau\in\Mf^\star(J)}  T(J,\sigma,\tau).
\end{align}
Note that the $\xi_\sigma^\star=0$ when $|J|>p$ since  $\Mf^\star(J)$ is empty in this case.

 \begin{lemma} \label{lemma_char_pol} 
 We have
\begin{align*}
\det(z- W)=  \sum_{k=0}^p \Big[\sum_{|J|=k,\sigma \in \Sym(J)} (-1)^{k+ \sign(\sigma)} \xi^\star_{\sigma} \Big]z^{n -k}.
\end{align*}
\end{lemma}

\proof  Using the definition of determinant, a direct computation gives
\begin{align*}
\det(z- W)=\det(z-M^*M)=  \sum_{k=0}^p \Big[\sum_{|J|=k,\sigma \in \Sym(J)} (-1)^{k+ \sign(\sigma)} \xi_{\sigma} \Big]z^{n -k}.
\end{align*}
Denote by $A_k$ the expression in the last brackets. 
According to (\ref{eq_xisigma}), $A_k$ is a polynomial in $t_{ij}$ and $\overline t_{ij}$.
Moreover, for each monomial in $A_k$,
 if some $j_0\in\llbracket 1,n\rrbracket$ appears as the second subscript of $t_{ij}$ or $\overline t_{ij}$, then it appears exactly one time for $t_{ij}$ and one time for $\overline t_{ij}$. We will call this property $\mathcal{P}_2$. The similar property for the first subscript will be called $\mathcal{P}_1$. 

In the same way, we obtain that $\det(z-MM^*)$ (here $M$ and $M^*$ were permuted) is a polynomial in $z$ whose coefficients satisfy $\mathcal{P}_1$. It follows from Lemma \ref{lemma_singvalue} that the polynomial $A_k$ satisfies both $\mathcal{P}_1$ and $\mathcal{P}_2$. Therefore, in the above definition of $A_k$, we can remove all monomials in $t_{ij}$ and $\overline t_{ij}$ which do not satisfy $\mathcal{P}_1$. This operation is equivalent to replacing $\xi_\sigma$ by $\xi_\sigma^\star$. The lemma follows.
\endproof

The following result gives us an upper bound for the expectation of the characteristic polynomial of $W$.

\begin{lemma} \label{lemma_det} 
We have 
$$\Ebf \big(\det(z- W)\big)=\sum_{J\subset \llbracket 1,n\rrbracket} (-1)^{|J|} |\Mf^\star(J)|z^{n-|J|}= (-1)^p p!   z^{n-p}  L^{(n-p)}_p(z).
$$
In particular, there is a universal  constant $c>0$ such that 
$$\big| \Ebf \big(\det(z- W)\big) \big|^2  \leq c   p^2 n!p!  z^{n-p}  e^{z}$$
for $z \in \R^+$ and  $n,p \in \N$ with $1\leq p\leq n$.
\end{lemma}
\proof 
We need the following claim whose proof will be given later.

\medskip\noindent
{\bf Claim.} Let $J$ be a subset of  $\llbracket1,n \rrbracket$ of cardinality $k \le p$ and $\sigma \in \Sym(J)$. 
Then $\Ebf (\xi^\star_{\sigma})= 0$ if $\sigma\not=\id_J$ and $\Ebf (\xi^\star_\sigma)=|\Mf^\star(J)|= p!/ (p-k)!$ if $\sigma=\id_J$.

\medskip

Using Lemma \ref{lemma_char_pol}, we have 
\begin{align*}
\Ebf (\det(z- W))= \sum_{k=0}^p \Big[\sum_{\sigma \in \Sym (J), |J|=k} (-1)^{k+ \sign(\sigma)} \Ebf (\xi^\star_{\sigma})\Big] z^{n -k}
\end{align*} 
which, by the above claim, is equal to 
\begin{eqnarray*}
\sum_{k=0}^p \Big[\sum_{|J|=k} (-1)^{k} \Ebf (\xi^\star_{\id_J})\Big] z^{n -k}
& = & \sum_{k=0}^p \Big[\sum_{|J|=k} (-1)^{k} |\Mf^\star(J)|\Big] z^{n -k} \\
& = & z^{n-p} \sum_{k=0}^p (-1)^k\binom{n}{k} \frac{p!}{(p-k)!} z^{p-k}.
\end{eqnarray*}
The first assertion in the lemma follows.

We apply Proposition \ref{prop_Laguerre} to $\alpha=n-p$ and to $z/n$ instead of $z$. We deduce that 
\begin{eqnarray*}
\big| \Ebf \big(\det(z- W)\big) \big|^2  & \leq & (p!)^2   z^{2n-2p}  |L^{(n-p)}_p(z)|^2 \\
& \lesssim &  (p!)^2   z^{2n-2p} \big[pn^{p\over 2} p^{-{p\over 2}} ({z/n})^{-{n-p\over 2}} e^{{z\over 2}} e^{-{n-p\over 2}}\big]^2 \\
& = &  p^2(p!)^2p^{-p}e^pn^ne^{-n} z^{n-p} e^z.
\end{eqnarray*}
Furthermore, by Stirling's formula, we have 
$$p!\simeq \sqrt{2\pi p} (p/e)^p \qquad \text{and} \qquad n!\simeq \sqrt{2\pi n} (n/e)^n\geq \sqrt{2\pi p} (n/e)^n.$$ 
The second assertion in the lemma follows easily.

It remains to prove the above claim. 
Denote by $\sigma_{1}, \cdots, \sigma_m$ the cycles of $\sigma.$ So we can write $\sigma= \sigma_1 \circ \cdots \circ \sigma_m.$ Note that here a fixed point of $\sigma$ is considered as a cycle of length $1.$ Observe that the $\xi^\star_{\sigma_{s}}$'s are mutually independent random variables. Thus $\Ebf (\xi^\star_{\sigma})= \prod_{1 \le s \le m} \Ebf (\xi^\star_{\sigma_s})$. 

\medskip\noindent
{\bf Case 1.} Assume that $\sigma \not = \id_J$. We want to show that $\Ebf(\xi^\star_\sigma)=0$.
 By (\ref{eq_xisigma}), it is enough to show that $\Ebf(T(J,\sigma,\tau))=0$.
 Since the $t_{ij}$'s are independent, using (\ref{eq_T}) we have 
 $$\Ebf(T(J,\sigma,\tau))=\Ebf\Big( \prod_{j\in J} \bar{t}_{\tau(j),  j} t_{\tau(\sigma^{-1}(j)),j}\Big)=\prod_{j \in J}  \Ebf  \big(\bar{t}_{\tau(j),  j} t_{\tau(\sigma^{-1}(j)),j}\big).$$
Since $\tau$ is injective and $\sigma$ is not the identity, there is a $j\in J$ such that $\tau(j)\not=\tau(\sigma^{-1}(j))$. For such a $j$, the last expectation satisfies
$$\Ebf  \big(\bar{t}_{\tau(j),  j} t_{\tau(\sigma^{-1}(j)),j}\big)= 
\Ebf  \big(\bar{t}_{\tau(j),  j} \big)\Ebf  \big(t_{\tau(\sigma^{-1}(j)),j}\big)=0$$
since  the $t_{ij}$'s are independent and have zero mean. The result follows.

\medskip\noindent
{\bf Case 2.} Assume now that $\sigma=\id_J$.  Using (\ref{eq_xisigma}) and arguing as above, we obtain that 
$$\Ebf (\xi^\star_{\sigma}) = \sum_{\tau\in\Mf^\star(J)}  \prod_{j \in J}  \Ebf  ( \bar{t}_{\tau(j), j} t_{\tau(j), j})=
 \sum_{\tau\in\Mf^\star(J)} 1 = |\Mf^\star(J)|$$
 since the variance of $t_{ij}$ is 1. Finally, recall that  $|J|=k$ and  $|\Mf^\star(J)|$  is the number of injective maps from $J$ to $\llbracket 1, p\rrbracket$. This number is equal to $p(p-1)\ldots (p-k+1)$. The lemma follows easily.  
\endproof   

We continue the proof of Proposition \ref{prop_det_W} and need to bound $\Ebf \big(|\det(z-W)|^2\big)$.
By Lemma \ref{lemma_char_pol}, we have 
\begin{align} \label{eq_det_sq}
\Ebf \big(|\det(z-W)|^2\big)= \sum_{J_1, J_2} \sum_{\sigma_1, \sigma_2} (-1)^{|J_1|+ |J_2|+ \sign(\sigma_1)+\sign(\sigma_2)} \Ebf(\xi^\star_{\sigma_1} \bar \xi^\star_{\sigma_2}) z^{2n -|J_1| -|J_2|},
\end{align}
where the sum is taken over the subsets $J_1, J_2$ of $\llbracket 1, n \rrbracket$ and 
the permutations $\sigma_1\in \Sym(J_1)$, $\sigma_2\in\Sym(J_2)$. Recall that we only need to consider the case where $|J_1|\leq p$ and $|J_2|\leq p$ since otherwise $\xi^\star_{\sigma_1}$ or $\xi^\star_{\sigma_2}$ vanishes. Note also that the proof of Lemma
\ref{lemma_det} suggests that fixed points of $\sigma_1$ and $\sigma_2$ may play an important role in our computation. This is the reason to use of the quantities introduced below.

Consider the set
\begin{align} \label{def_Mf-double}
\Mf^{\star}(J_1,J_2):= \Big\{(\tau_1,\tau_2)\in \Mf^\star(J_1)\times \Mf^\star(J_2) : \  \tau_1(j)\not=\tau_2(j) \text{ for every } j\in J_1\cap J_2 \Big\}
\end{align}
and define
 \begin{align} \label{def_R1}
 R^{[1]}(n,p,z):= \sum_{J_1, J_2} (-1)^{|J_1|+ |J_2|}|\Mf^{\star}(J_1,J_2)|  z^{2n - |J_1| - |J_2|},
 \end{align}
  where the sum is  taken over the subsets $J_1$ and $J_2$ of  $\llbracket 1, n\rrbracket$.
Note that $|\Mf^{\star}(J_1,J_2)|$ only depends on $|J_1|, |J_2|,p$ and it vanishes when $|J_1|>p$ or $|J_2|>p$. 
 Define also
 \begin{align} \label{def_R2}
 R^{[2]}(W,z):= \sum_{J_1, J_2} (-1)^{|J_1|+ |J_2|} \Ebf\big(\xi^\star_{\id_{J_1}} \bar\xi^\star_{\id_{J_2}}\big)  z^{2n - |J_1| - |J_2|}.
 \end{align}

\begin{lemma} \label{lemma_R2}
We have
$$R^{[2]}(W,z)= \sum_{J \subset \llbracket 1, n\rrbracket} R^{[1]}(n- |J|,p- |J|,z)  \Big[\sum_{\tau\in\Mf^\star(J)}\prod_{j \in J} \Ebf\big( |t_{\tau(j), j}|^4\big)\Big].$$
In particular, we have 
$$|R^{[2]}(W,z)| \le  e^{\beta} \,  n! \, p! \sum_{k=0}^p   \frac{|R^{[1]}(n- k, p- k,z)|}{(n-k)! \, (p-k)!} \cdot     $$
\end{lemma}
\proof We prove the first assertion. 
By (\ref{eq_T}) and (\ref{eq_xisigma}),  we have 
\begin{align} \label{eq_xixi}
\xi^\star_{\id_{J_1}} \bar \xi^\star_{\id_{J_2}}= \sum_{\tau_1\in \Mf^\star(J_1)}\sum_{\tau_2\in\Mf^\star(J_2)} \prod_{j \in J_1} |t_{\tau_1(j),j}|^2  \, \prod_{j \in J_2}  |t_{\tau_2(j), j}|^2.
\end{align}
Consider the summand in the right hand side of this identity. Observe that term  $ |t_{i j}|^4$ only appears when $j \in J_1 \cap J_2$ and  $\tau_1(j)= \tau_2(j).$ 
Denote by  $J$ the set of $j\in J_1 \cap J_2$ such that $\tau_1(j)=\tau_2(j)$ and $\tau$ is the restriction of $\tau_1$ and $\tau_2$ to $J$.  Define also $J_1':=J_1\setminus J$, $J_2':=J_2\setminus J$ and $\tau_1',\tau_2'$ the restrictions of $\tau_1,\tau_2$ to $J_1',J_2'$ respectively.
So the summand in the right hand side of   (\ref{eq_xixi})  is equal to 
$$\prod_{j \in J}|t_{\tau(j), j}|^4 \prod_{j \in J_1'} |t_{\tau'_1(j), j}|^2 \prod_{j \in J_2'} |t_{\tau'_2(j), j}|^2.$$

Note that $J_1'$ and $J_2'$ are subsets of $\llbracket 1,n\rrbracket\setminus J$ and the later set is of cardinality $n':=n-|J|$. Moreover, $\tau_1'$ and $\tau_2'$ have images in the set $\llbracket 1, p\rrbracket\setminus \tau(J)$ which is of cardinality $p':=p-|J|$ and they are not equal at any point of $J_1'\cap J_2'$. Therefore, for each pair $(i,j)$ with $j\in J_1'\cup J_2'$, if $|t_{ij}|^2$ appears in the last product, it appears exactly 1 time. Thus, the expectation of this product is equal to 
$$e(J,\tau):=  \prod_{j \in J}\Ebf\big( |t_{\tau(j), j}|^4\big).$$

For the right hand side of  (\ref{def_R2}), observe that
$$(-1)^{|J_1|+ |J_2|}  z^{2n - |J_1| - |J_2|} = 
(-1)^{|J_1'|+ |J_2'|} z^{2n' - |J_1'| - |J_2'|}.$$
Furthermore, if we fix $J_1, J_2, J$ and $\tau\in \Mf^\star(J)$, then the family of all pairs $\tau_1',\tau_2'$ is similar to the set $\Mf^\star(J_1,J_2)$ defined in (\ref{def_Mf-double}). 
We denote it by $\widetilde\Mf^\star(J_1',J_2')$.
The main difference is that $J'_1, J'_2$ are subsets of $\llbracket 1, n\rrbracket \setminus J$ instead of $\llbracket 1, n\rrbracket$ and $\tau'_1,\tau_2'$ have images in $\llbracket 1, p\rrbracket \setminus\tau(J)$ instead of $\llbracket 1, p\rrbracket $. This doesn't cause any difficulty in our computation because we will only use the cardinality of this set.
Using (\ref{def_R1}) for $n',p'$ instead of $n,p$, we obtain that
\begin{eqnarray*}
R^{[2]}(W,z) & = &  \sum_{J \subset \llbracket 1, n\rrbracket}  \sum_{\tau\in\Mf^\star(J)} e(J,\tau) \sum_{J'_1,J'_2} |\widetilde\Mf^\star(J_1',J_2')| (-1)^{|J_1'|+ |J_2'|} z^{2n' - |J_1'| - |J_2'|} \\
&  = & 
\sum_{J \subset \llbracket 1, n\rrbracket}  \sum_{\tau\in\Mf^\star(J)} e(J,\tau) R^{[1]}(n',p',z).
\end{eqnarray*}
The first assertion in the lemma follows easily.

We now prove the second assertion using the first one. Observe that the number of subsets $J$ of $\llbracket 1,n\rrbracket$ of cardinality $k$ is $\binom{n}{k}$ and for such a set, the cardinality of $\Mf^\star(J)$ is $p!/(p-k)!$. Therefore,  
using the first assertion and the fact that $\Ebf(|t_{ij}|^4)\leq \beta$, we get
\begin{eqnarray*}
|R^{[2]}(W,z)|  & \le &   \sum_{k=0}^p \binom{n}{k} \frac{p!}{(p-k)!} \beta^k  |R^{[1]}(n- k, p- k,z)| \\
& = & n! \, p! \sum_{k=0}^p  \frac{\beta^k}{k!} \,  \frac{|R^{[1]}(n- k, p- k,z)|}{(n-k)! \, (p-k)!} \\ 
&\leq & e^\beta n! \, p! \sum_{k=0}^p  \frac{|R^{[1]}(n- k, p- k,z)|}{(n-k)! \, (p-k)!} \cdot
\end{eqnarray*}
This ends the proof of the lemma.
\endproof

In order to bound $\Ebf (|\det(z-W)|^2)$, we will need to bound  $R^{[1]}(n',p', z)$ for $0 \le n' \le n$ and $0 \le p' \le n'.$ The following lemma is a crucial step.

\begin{lemma} \label{lemma_W_R1}
We have
$$\big|\Ebf \big(\det(z-W)\big)\big|^2 = n! \, p! \sum_{k=0}^p \frac{ R^{[1]}(n-k, p-k,z)}{k!\, (n-k)! \, (p-k)!} \cdot$$
\end{lemma} 
 \proof 
 We have by Lemma \ref{lemma_det} that
 $$\Ebf \big(\det(z- W)\big)= \sum_J (-1)^{|J|} |\Mf^\star(J)| z^{n -|J|}$$
 which implies that
 \begin{align} \label{eq_W_MM}
 |\Ebf \det(z- W)|^2= \sum_{J_1, J_2}   (-1)^{|J_1|+ |J_2|}|\Mf^\star(J_1)\times \Mf^\star(J_2)|  z^{2n -|J_1| -|J_2|}.
 \end{align}
In order to relate the last sum to $R^{[1]}(n,p,z)$, we will compute  $|\Mf^\star(J_1)\times \Mf^\star(J_2)|$ in a suitable way.

Consider $(\tau_1,\tau_2)\in \Mf^\star(J_1)\times \Mf^\star(J_2)$. Using the notations introduced in 
 the proof of Lemma \ref{lemma_R2}, we have 
 $$|\Mf^\star(J_1)\times \Mf^\star(J_2)|=\sum_{J,\tau} \sum_{J_1',J_2'} |\widetilde\Mf^\star(J_1',J_2')|.$$
For the right hand side of 
(\ref{eq_W_MM}), we also observe that
$$(-1)^{|J_1|+ |J_2|} z^{2n -|J_1| -|J_2|} = (-1)^{|J_1'|+ |J_2'|} z^{2n' -|J_1'| -|J_2'|}.$$
This arguments, together with (\ref{def_R1}), imply 
$$ |\Ebf \det(z- W)|^2  =   \sum_{J,\tau} \sum_{J_1',J_2'}  (-1)^{|J_1'|+ |J_2'|}|\widetilde\Mf^\star(J_1',J_2')|  z^{2n' -|J_1'| -|J_2'|}  =   \sum_{J,\tau} R^{[1]}(n',p',z).$$
Finally, recall that the number of subsets $J$ of $\llbracket 1,n\rrbracket$ of cardinality $k$ is $\binom{n}{k}$ and for such a set $J$ the number of $\tau$ in $\Mf^\star(J)$ is $p!/(p-k)!$. We easily deduce the lemma from the last identities.
 \endproof

 \begin{lemma} \label{lemma_R_bound} 
There is a universal constant $c>0$ such that
$$ |R^{[1]}(n,p,z)|    \leq c(p+1)^3  n!  p!  z^{n-p} e^z \quad \text{and} \quad  |R^{[2]}(W,z)|   
 \leq c (p+1)^4 e^{\beta}   n!  p!  z^{n-p} e^z$$
for  $z \in \R^+$,  $n\geq 1$ and $0\leq p \leq n$.
 \end{lemma}

\proof 
The second estimate is a direct consequence of the first one and Lemma \ref{lemma_R2}. 
We will prove the first estimate by induction in $p$. Therefore, we will write $W_{p,n}$ instead of $W$ in order to mention the size of the matrix $M$ which is used in the definition of $W$. Note that $R^{[1]}(n,0,z)=z^{2n} \le n! z^n  e^z.$  
So the estimate holds for $p=0$.
By Lemma \ref{lemma_W_R1} for $p=1$, we have 
$$|R^{[1]}(n,1,z)|= |\Ebf \det(z-W_{1,n})|^2-  n R^{[1]}(n-1,0,z) \le  |\Ebf \det(z-W_{1,n})|^2.$$ 
Thus, Lemma \ref{lemma_det} implies the desired estimate for $p=1.$

Now, let $p\geq 2$ be an integer and assume that the desired estimate  holds for every $p'$ such that $0\leq p' \le p-1$ and for a universal constant $c>0$ large enough. We need to prove this estimate for $p$.  Define
$$r_{n,p}(z):=\frac{ R^{[1]}(n,p,z)}{n!  p! z^{n-p}   e^z} \qquad \text{and} \qquad e_{n,p}(z):= \frac{ |\Ebf \det(z-W_{p,n})|^2 }{n! p!  z^{n -p} e^z} \cdot$$
By the induction hypothesis,  we have for $1 \le k \le p$
\begin{align} \label{ineq_rnp}
|r_{n-k,p-k}(z)| \leq c(p-k+1)^3\leq cp^3.
\end{align}
We need to check that $|r_{n,p}(z)|\leq c(p+1)^3$. Note also that by Lemma \ref{lemma_det} we have for $c$ large enough
\begin{align} \label{ineq_enp}
 0\leq e_{n,p}(z) \leq c p^2 \qquad \text{and} \qquad 0\leq e_{n-1,p-1}(z) \leq c p^2.
\end{align}

By Lemma \ref{lemma_W_R1}, we have
\begin{align} \label{eq_r_np}
r_{n,p}(z)= e_{n,p}(z)- \sum_{k=1}^p \frac{r_{n-k, p-k}(z)}{k!}=e_{n,p}(z)- r_{n-1,p-1}(z)-\sum_{k=2}^p \frac{r_{n-k, p-k}(z)}{k!}\cdot 
\end{align}
The same lemma applied to $n-1,p-1$ instead of $n,p$ gives
$$r_{n-1,p-1}(z)= e_{n-1,p-1}(z)- \sum_{k=1}^{p-1} \frac{r_{n-1-k, p-1-k}(z)}{k!}= e_{n-1,p-1}(z)- \sum_{k'=2}^{p} \frac{r_{n-k', p-k'}(z)}{(k'-1)!}\cdot $$
Substituting the obtained value of $r_{n-1,p-1}(z)$ to (\ref{eq_r_np}) gives
$$r_{n,p}(z)= e_{n,p}(z) - e_{n-1,p-1}(z)+ \sum_{k=2}^p \Big[\frac{1}{(k-1)!}- \frac{1}{k!}\Big]  r_{n-k, p-k}(z). $$
Finally, combining this with (\ref{ineq_rnp}) and (\ref{ineq_enp}), we obtain
$$|r_{n,p}(z)| \leq c p^2 +  \sum_{k=2}^p  \Big[\frac{1}{(k-1)!}- \frac{1}{k!}\Big] cp^3 \leq cp^2+cp^3\leq c(p+1)^3.$$
This ends the proof of the lemma.
\endproof

Recall that we have to bound the expectation of $|\det (z-W)|^2$. By (\ref{eq_det_sq}), we need to study
 $\Ebf (\xi^\star_{\sigma_1} \bar \xi^\star_{\sigma_2})$ for general $\sigma_1$ and $\sigma_2.$ 
 So we need to consider $T(J,\sigma,\tau)$ defined in (\ref{eq_T}) for $J\subset \llbracket 1,n\rrbracket$ of cardinality at most $p$, $\sigma\in\Sym(J)$ and $\tau\in\Mf^\star(J)$.
 In order to study this quantity, 
 we will use a graph associated to $J,\sigma$ and $\tau$ as in \cite[section 3.1.2]{BS}.

Draw two parallel lines $L^+$ and $L^-$  in $\R^2$. Consider the points of abscissas $1,2,\ldots, p$ in the upper line $L^+$, corresponding to the index $i$ of $t_{ij}$ and the points of abscissas $1,2,\ldots, n$ in the lower line $L^-$, corresponding to the index $j$ of $t_{ij}$. 
These points will be the vertices of the following graphs.
Define the oriented graph $\Gamma(J,\sigma,\tau)$  by 

\smallskip
- joining the point of abscissa  $\tau(j)$ in $L^+$
with the point of abscissa $j\in J$ in $L^-$  by a down edge, and  

\smallskip
- joining the point of abscissa  $j\in J$ in $L^-$ with the point of abscissa
$\tau(\sigma^{-1}(j))$ in $L^+$ by  an up edge. 

\smallskip\noindent
So each vertex of this graph is of degree 2 and belongs to exactly one up edge and one down edge (we don't consider the points in $L^+\cup L^-$ which don't belong to any edge).

Observe that if $\sigma'$ is a cycle of $\sigma$, $J'$ its support and $\tau'$ the restriction of $\tau$ to $J'$, then $\Gamma(J',\sigma',\tau')$ is a cycle and also a connected component of $\Gamma(J,\sigma,\tau)$.  In particular, if $j$ is a fixed point of $\sigma$ then the point of abscissa $j$ in $L^-$, the point of abscissa $\tau(j)$ in $L^+$ and the up and down edges joining them constitute a connected component of $\Gamma(J,\sigma,\tau)$. We call them components of fixed points. They are the only components in which two vertices are joined by an up edge and a down edge. We also see that there is a 1:1 correspondence between $(J,\sigma,\tau)$ and $\Gamma(J,\sigma,\tau)$. 

Observe also that $\bar t_{ij}$ appears in the monomial $T(J,\sigma,\tau)$ if and only if there is a down edge of  $\Gamma(J,\sigma,\tau)$ joining the point of abscissa $i$ in $L^+$ to the point of abscissa $j$ in $L^-$. Similarly, $t_{ij}$ appears in this monomial  if and only if there is a up edge $\Gamma(J,\sigma,\tau)$ joining 
the point of abscissa $j$ in $L^-$ to the point of abscissa $i$ in $L^+$. 
 The following lemma is important for our computation.

 \begin{lemma} \label{lemma_TTbar} 
 Let $J_1$ and $J_2$ be subsets of length at most $p$ of  $\llbracket 1,n\rrbracket$. Let $\sigma_k \in \Sym (J_k)$ and
 $\tau_k \in \Mf^\star (J_k)$ for $k=1,2$. 
 We say that $(J_1,\sigma_1,\tau_1)$ and $(J_2,\sigma_2,\tau_2)$ are compatible if every component of non-fixed points in $\Gamma(J_1,\sigma_1,\tau_1)$
is a component of $\Gamma(J_2,\sigma_2,\tau_2)$ regardless the orientation, and vice versa.
 Then if $(J_1,\sigma_1,\tau_1)$ and $(J_2,\sigma_2,\tau_2)$ are not compatible, 
 we have
 $$\Ebf\big(T(J_1,\sigma_1,\tau_1)\overline{T(J_2,\sigma_2,\tau_2)}\big)=0.$$
 In particular, we have 
$\Ebf (\xi^\star_{\sigma_1} \bar\xi^\star_{\sigma_2}) =0$ in the following cases :
 \begin{enumerate}
 \item $\sigma_1$ contains a cycle $l$ of length at least $2$ such that neither $l$ nor $l^{-1}$ is a cycle of $\sigma_2$;
 \item $\sigma_2$ contains a cycle $l$ of length at least $2$ such that neither $l$ nor $l^{-1}$ is a cycle of $\sigma_1$;
 \end{enumerate}
 \end{lemma}

 \proof 
Observe that in the situation of (1) or (2), $(J_1,\sigma_1,\tau_1)$ and $(J_2,\sigma_2,\tau_2)$ are not compatible for all $\tau_1$ and $\tau_2$.
By (\ref{eq_xisigma}), the second assertion is clearly a consequence of the first one. We prove now the first assertion.
Recall from (\ref{eq_T}) that 
$$T(J_1,\sigma_1,\tau_1)\overline{T(J_2,\sigma_2,\tau_2)}=\prod_{j \in J_1} \bar t_{\tau_1(j),j} t_{\tau_1(\sigma_1^{-1}(j)), j}\prod_{j \in J_2} t_{\tau_2(j),j} \bar t_{\tau_2(\sigma_2^{-1}(j)), j}.$$
This is a monomial in $t_{ij}$'s and $\bar t_{ij}$'s. 
By hypothesis, the two graphs are not compatible, there are two vertices which are joined by exactly one edge of the union of these two graphs.  So for some $i,j$, the  the total degree of this monomial in $t_{ij}$ and $\bar t_{ij}$ is 1.
Since the $t_{ij}$'s are independent and of zero mean, we deduce that the expectation of the considered monomial is also 0. The lemma follows.
\endproof

According to Lemma \ref{lemma_TTbar}, in our computation of expectations, we only need to consider compatible $(J_1,\sigma_1,\tau_1)$ and $(J_2,\sigma_2,\tau_2)$. Let $J_k^F$  denote the sets of fixed points of $\sigma_k$. We have the following properties :

\smallskip
(P1) We have $J_1\setminus J_1^F=J_2\setminus J_2^F$. Denote by $J_0$ this set which is contained in $J_1\cap J_2$. Denote by $\sigma_{0k},\tau_{0k}$ the restrictions of $\sigma_k,\tau_k$ to $J_0$. 

\smallskip
(P2) If $l$ is a cycle in $\sigma_{01}$ then either $l$ or $l^{-1}$ is a cycle in $\sigma_{02}$, and vice versa. Note that the length of $l$ is at least 2 or equivalently $\sigma_{01}$ and $\sigma_{02}$ have no fixed point. So  $\tau_{02}$ is uniquely determined by 
$\tau_{01},\sigma_{01}$ and $\sigma_{02}$.

\smallskip
(P3) Then the graphs $\Gamma(J_0,\sigma_{01},\tau_{01})$ and 
$\Gamma(J_0,\sigma_{02},\tau_{02})$ have the same support, that is, they are equal if we don't consider the orientation. In particular, we have $\tau_{01}(J_0)=\tau_{02}(J_0)$. If $I_0$ denotes the last set, then 
$\tau_k(J^F_k)$ is contained  in $\llbracket 1,n\rrbracket\setminus I_0$.

\smallskip

Denote by $M[J_0,\tau_{01},\tau_{02}]$ the matrix obtained from $M$ by removing the columns of index $j\in J_0$ and the lines of index $i\in I_0$. This is a matrix of size $p_0\times n_0$ with $p_0:=p-|J_0|$ and $n_0:=n-|J_0|$.  Denote also by $W[J_0,\tau_{01},\tau_{02}]$ the covariance matrix associated with $M[J_0,\tau_{01},\tau_{02}]$.

\begin{lemma} \label{lemma_det_sq_R2}
With the above notation, we have that $\Ebf (|\det(z-W)|^2)$ is equal to
$$\sum_{J_0, \sigma_{01},\sigma_{02},\tau_{01},\tau_{02}} (-1)^{\sign(\sigma_{01})+\sign(\sigma_{02})} \Ebf\big(T(J_0,\sigma_{01},\tau_{01})\overline{T(J_0,\sigma_{02}, \tau_{02})}\big) R^{[2]}(W[J_0,\tau_{01},\tau_{02}],z),$$
where $J_0,\sigma_{01},\sigma_{02}, \tau_{01},\tau_{02}$ satisfy  (P1), (P2), (P3) and $R^{[2]}(W[J_0,\tau_{01},\tau_{02}],z)$ is defined as in (\ref{def_R2}) but for $W[J_0,\tau_{01},\tau_{02}]$ instead of $W$.
\end{lemma}
\proof
Denote by $\tau_k^F$ the restriction of $\tau_k$ to $J_k^F$.
Using the above discussion, the identities (\ref{eq_det_sq}), (\ref{eq_xixi}) and the fact that the $t_{ij}$'s are independent,  
we see that $\Ebf (|\det(z-W)|^2)$ is equal to
\begin{align*}
\sum_{J_0, \sigma_{01},\sigma_{02},\tau_{01},\tau_{02}}  \sum_{J_1^F, J_2^F,\tau_1^F,\tau_2^F} (-1)^{|J_1|+|J_2|+\sign(\sigma_1)+\sign(\sigma_2)}\Ebf\big(T(J_0,\sigma_{01},\tau_{01})\overline{T(J_0,\sigma_{02}, \tau_{02})}\big) \\
\Ebf\big(T(J_1^F,\id_{J_1^F},\tau_1^F) \overline{T(J_2^F,\id_{J_2^F},\tau_2^F)}\big)z^{2n-|J_1|-|J_2|}.
\end{align*}
Observe that 
$$(-1)^{|J_1|+|J_2|+\sign(\sigma_1)+\sign(\sigma_2)}z^{2n-|J_1|-|J_2|}=(-1)^{|J_1^F|+|J_2^F|} (-1)^{\sign(\sigma_{01})+\sign(\sigma_{02})}z^{2n_0-|J_1^F|-|J_2^F|}.$$
Thus, $\Ebf (|\det(z-W)|^2)$ is equal to
\begin{align*}
\sum_{J_0, \sigma_{01},\sigma_{02},\tau_{01},\tau_{02}}  (-1)^{\sign(\sigma_{01})+\sign(\sigma_{02})} \Ebf\big(T(J_0,\sigma_{01},\tau_{01})\overline{T(J_0,\sigma_{02}, \tau_{02})}\big) \hspace{4.5cm}\  \\ 
\sum_{J_1^F, J_2^F,\tau_1^F,\tau_2^F} (-1)^{|J_1^F|+|J_2^F|}
\Ebf\big(T(J_1^F,\id_{J_1^F},\tau_1^F) \overline{T(J_2^F,\id_{J_2^F},\tau_2^F)}\big)z^{2n_0-|J_1^F|-|J_2^F|}.
\end{align*}
It follows from (\ref{def_R2}) and (\ref{eq_xisigma}), applied to $W[J_0,\tau_{01},\tau_{02}]$ instead of $W$, that the last summation is equal to $R^{[2]}(W[J_0,\tau_{01},\tau_{02}],z)$. This implies the lemma. 
\endproof

\noindent
{\bf End of the proof of Proposition \ref{prop_det_W}.}
We will apply Lemmas \ref{lemma_R_bound} and \ref{lemma_det_sq_R2}. First, observe that since $t_{ij}$ is of variance 1, we have $|\Ebf(t_{ij}^2)|=|\Ebf(\overline t_{ij}^2)|\leq 1$ and $\Ebf(|t_{ij}|^2)=1$. Therefore, since for each pair $(i,j)$, the total degree  in $t_{ij}$ and $\bar t_{ij}$ of the monomial 
$T(J_0,\sigma_{01},\tau_{01})\overline{T(J_0,\sigma_{02}, \tau_{02})}$  is  0 or 2, 
we deduce that
$$|\Ebf \big(T(J_0,\sigma_{01},\tau_{01})\overline{T(J_0,\sigma_{02}, \tau_{02})}\big)|\leq 1.$$
By Lemmas  \ref{lemma_R_bound} and  \ref{lemma_det_sq_R2}, we have 
\begin{eqnarray*}
\Ebf (|\det(z-W)|^2) &\lesssim & \sum_{J_0, \sigma_{01},\sigma_{02},\tau_{01},\tau_{02}}|R^{[2]}(W[J_0,\tau_{01},\tau_{02}],z)| \\
& \lesssim & \sum_{J_0, \sigma_{01},\sigma_{02},\tau_{01},\tau_{02}} (p-|J_0|+1)^4 e^\beta (n-|J_0|)! (p-|J_0|)! z^{n-p}e^z \\
& \lesssim & p^4e^\beta  z^{n-p}e^z\sum_{J_0, \sigma_{01},\sigma_{02},\tau_{01},\tau_{02}}  (n-|J_0|)! (p-|J_0|)! .
\end{eqnarray*}

The number of subsets $J_0$ of $\llbracket 1,n\rrbracket$ of cardinality $k$ is $\binom{n}{k}$. When such a $J_0$ is fixed, there are $p!/(p-k)!$ choices for the injective map $\tau_{01}$ from $J_0$ to $\llbracket 1, p\rrbracket$. 
By (P3), the map $\tau_{02}$ is uniquely determined by $\tau_{01}, \sigma_{01}$ and $\sigma_{02}$. 
We also have the following property that we will prove later.

\medskip\noindent
{\bf Claim.} When the set $J_0$ of cardinality $k$ is fixed, there are less than $(k+1)!$ choices for the pair $(\sigma_{01},\sigma_{02})$ satisfying the above property (P2). 

\medskip\noindent
This, together with  the Stirling's formula $n!\simeq \sqrt{2n\pi}({n\over e})^n$ and $p!\simeq \sqrt{2p\pi}({p\over e})^p$, imply that 
\begin{eqnarray*}
\Ebf (|\det(z-W)|^2) &\lesssim & p^4e^\beta  z^{n-p}e^z\sum_{k=0}^p \binom{n}{k}{p!\over (p-k)!} (n-k)!(p-k)! (k+1)!\\
& \lesssim & p^6n! p!e^\beta  z^{n-p}e^z \ \lesssim \ p^{13/2}n^{1/2} e^\beta n^n p^p z^{n-p}e^{z-n-p} .
\end{eqnarray*}
Hence, 
$$\Ebf (|\det(z-n^{-1}W)|^2) =n^{-2n} \Ebf (|\det(nz-W)|^2)\lesssim p^{13/2}n^{1/2} e^\beta n^{-p} p^p z^{n-p}e^{nz-n-p} .$$
A direct computation using Proposition \ref{prop_potential_MP} shows that the last expression is equal to $p^{13/2}n^{1/2} e^\beta e^{2nu_\phi(z)}$. This gives the desired estimate of the proposition.

In what follows, we prove the above claim. We will show that there are at most $(k+1)!$ choices for $(\sigma_{01},\sigma_{02})$ satisfying the first part of (P2), that is, we accept fixed points. 
Let $\mathfrak{s}(k,m)$ be the unsigned Stirling's number which is also the number of permutations of $J_0$ having exactly $m$ cycles \cite[p.234]{Comtet}. We have  for $x \in \R$
\begin{align*} 
\sum_{m=1}^k \mathfrak{s}(k,m) x^m= x(x+1) \ldots (x+k-1),
\end{align*}
see \cite[Th. A, p.213]{Comtet}. 

Observe that when we fix a permutation $\sigma_{01}$ of $J_0$ which has exactly $m$ cycles (a fixed point is also considered as a cycle), there are at most $2^m$ choices for $\sigma_{02}$ satisfying the first part of (P2). Therefore, the total number of choices  for $(\sigma_{01},\sigma_{02})$ is at most 
\begin{align*}
 \sum_{m=1}^k \mathfrak{s}(k,m) 2^m= 2(2+1) \ldots (k+1)= (k+1)!
\end{align*}
The proof of the proposition is now complete.
\hfill $\square$

\medskip

We continue the proof of Theorem \ref{th_main}. Let $L$ be a compact interval in $\R$. Let $K$ be a compact interval contained in 
$L$.  Let $\mu_0$ be a given probability measure whose support is contained in $K$ and 
$u_0$ its logarithmic potential. 
We assume that there is a constant $\kappa\geq 1$ such that for $x_1,x_2\in K$ we have 
$$|u_0(x_1)-u_0(x_2)|\leq \kappa  |x_1-x_2|.$$

Let $\Pc_0^p$ denote the set of all monic complex polynomials of one variable and degree $p$. Consider  a probability measure on $\Pc_0^p$ with $p\geq 1$. 
If $Q$ is a polynomial of degree $p$ and $z_1,\ldots, z_p$ are its zeros, define 
$$\mu_Q:={1\over p} \sum_{k=1}^p \delta_{z_k}.$$
This is the probability measure equidistributed on the zero set of $Q$. 
We recall the following large deviation theorem from \cite{Dinh} that we will use later for $L:=[0,4]$, $K:=[a,b]$ and $\mu_0:=\mu_\phi^+$.

\begin{theorem} \label{th_dinh}
Let $L, K,\mu_0,u_0$ and $\kappa$ be as above. 
Consider a probability measure on the set $\Pc^p_0$ of all monic polynomials of degree $p$ and a constant $c_p\geq 1$  
such that the expectation of $|Q(z)|^2$   for $Q\in \Pc^p_0$ satisfies
$$\Ebf(|Q(z)|^2)\leq  c_p e^{2 p u_0(z)} \quad \text{for all} \quad  z\in K.$$ 
Then there are positive constants $A_1$ and $A_2$ depending only on $L$ such that
$$\Prob\Big\{Q\in \Pc^p_0, \quad \dist(\mu_Q,\mu_0) \geq \delta\Big\} \leq A_1 p^2 \kappa c_p e^{-A_2\delta p}.$$
\end{theorem}
\proof
The only difference between this statement and a main result of  \cite{Dinh} is the precision on the constants. We will sketch here the proof in order to see this point. Choose positive constants $A$ large enough and $A_2$ small enough which only depend on $L$. 
Let $\Sigma$ be a set of about $2A\kappa p^2$ points equidistributed on $K$. Let $\Ec$ be a set of $Q$ such that 
$|Q(a)|^2 \geq   A^{-2}e^{2 p u_0(a)} e^{A_2\delta p}$ for some $a\in \Sigma$. It is not difficult to see that the probability that a polynomial $Q$ belongs to  $\Ec$ is at most $2A^3 p^2 \kappa c_p e^{-A_2\delta p}$. It is enough to check for $Q\not\in\Ec$ that 
$\dist(\mu_Q,\mu_0) \leq \delta$. By Proposition \ref{prop_dist_sup}, we only need to check that $|Q(z)|\leq e^{pu_0(z)}e^{{1\over 2} A_2 p\delta}$ for $z\in K$ since $A_2$ is small enough. Thus, we only have to prove this inequality for $b\in K$ such that $|Q(z)|e^{-pu_0(z)}$ is maximal when $z=b$. Recall that we already have $|Q(a)| \leq   A^{-1}e^{p u_0(a)} e^{{1\over 2}A_2\delta p}$ for all $a\in \Sigma$.

If $m$ denotes the length of $K$, there is a point $a\in\Sigma$ such that $|a-b|\leq m A^{-1}\kappa^{-1}p^{-2}$. Choose an interval $K'\subset K$ of length $m':=\min(m,1)\kappa^{-1}p^{-1}$ containing $a$ and $b$ (this is possible when we choose $A$ larger than 1 and  larger than the length of $L$).  By definition of $\kappa$ and $b$ we get
$$\max_{K'} pu_0-\min_{K'} pu_0\leq \min(m,1) \quad \text{and hence} \quad \max_{K'} |Q| \leq e|Q(b)|.$$
Using the classical Markov brother's theorem and a simple change of coordinate, we have
$$\max_{K'}|Q'| \leq m'^{-1} p \max_{K'} |Q| \leq e\min(m,1)^{-1}\kappa p^2 |Q(b)|.$$ 
It follows that 
$$|Q(a)|\geq |Q(b)|-|Q(a)-Q(b)| \geq |Q(b)|- e\min(m,1)^{-1}\kappa p^2 |Q(b)||a-b| \geq {1\over 2} |Q(b)|$$
since $A$ is large. Thus,
$$|Q(b)|e^{-pu_0(b)} \leq 2  |Q(a)|e^{-pu_0(a)+1}  \leq 2e A^{-1} e^{{1\over 2}A_2\delta p}\leq e^{{1\over 2}A_2\delta p}.$$
We used in the first inequality that $pu_0(a)- pu_0(b) \leq 1$. The result follows.
\endproof

\noindent
{\bf End of the proof of Theorem \ref{th_main}.}
By Lemma \ref{lemma_singvalue}, the matrix $W$ has 0 as an eigenvalue of multiplicity at least $n-p$. So to each $n^{-1}W=n^{-1} M^*M$ we can associate the polynomial $Q(z):=z^{-n+p} \det (z-n^{-1}W)$ which is an element of $\Pc_0^p$.
This induces a probability measure on $\Pc_0^p$ because  the entries of $M$ are random variables. Moreover, by Proposition \ref{prop_det_W} and the identity before Proposition \ref{prop_potential_MP}, we have for $z\in[a,b]$
$$\Ebf(|Q(z)|^2) \leq cp^{13/2} n^{1/2}  e^\beta e^{2pu_\phi^+(z)}\leq cn^7 e^\beta e^{2pu_\phi^+(z)}.$$

Denote by $\mu_{p,n}^+$ the probability measure $\mu_Q$ with $Q$ associated to $W$ as above. 
We apply 
Theorem \ref{th_dinh} for $L:=[0,4]$, $K:=[a,b]$,  $\mu_0:=\mu_\phi^+$, $u_0:=u_\phi^+$ and for $\gamma\phi^{-1}\delta$ instead of $\delta$, where $\gamma>0$ is a fixed constant small enough. 
It suffice to choose $\kappa$ larger than the sup-norm of the derivative of $u^+_\phi$ on $[a,b]$. So
by Proposition \ref{prop_potential_MP}, if $\phi=1$ we can choose $\kappa$ any constant larger than 1. Otherwise, since $a\geq (1-\phi)^2/4$, we can choose any $\kappa$ larger than a constant times $\phi^{-1}(1-\phi)^{-1}$. Observe that in the last case we have $\phi^{-1}(1-\phi)^{-1}\leq n^2$ because both $\phi$ and $1-\phi$ are larger than $1/n$. 
So in any case, we can choose $\kappa$ equal to a big enough constant times $n^2$. 

Let $A_1$ and $A_2$ be as in  Theorem \ref{th_dinh}. Take $A_2':=A_2\gamma$ and choose $A_1'$ equal to a large constant times $A_1$. In the present context, Theorem \ref{th_dinh} implies the existence of 
  a set $\Ec_{p,n}(\delta)$ of $(p\times n)$-matrices such that 
$$\Prob(M\in \Ec_{p,n}(\delta))\leq 
A_1'n^{11} e^\beta e^{-A_2p(\gamma\phi^{-1}\delta)} = 
A_1'n^{11} e^\beta e^{-A_2'\delta n}$$
and also 
$$\dist(\mu_{p,n}^+,\mu_\phi^+)\leq \gamma\phi^{-1}\delta \quad \text{for} \quad M\not\in \Ec_{p,n}(\delta).$$
Since $\mu_{p,n}=(1-\phi)\delta_0+\phi\mu_{p,n}^+$ and $\mu_\phi=(1-\phi)\delta_0+\phi\mu_\phi^+$, the last inequality implies that
$$\dist(\mu_{p,n},\mu_\phi)\leq \gamma\delta\leq\delta  \quad \text{for} \quad M\not\in \Ec_{p,n}(\delta).$$
This is one of the desired inequalities in the theorem.

Now, we can see in the proof of Theorem \ref{th_dinh} (applied to $\gamma\phi^{-1}\delta$ instead of $\delta$) that one can choose $\Ec_{p,n}(\delta)$ so that if $M\not\in \Ec_{p,n}(\delta)$ then
$$|Q(z)|\leq e^{pu_\phi^+(z)} e^{{1\over 2} A_2p(\gamma\phi^{-1}\delta)} \qquad \text{for}\quad z\in K$$
or equivalently
$${1\over p}\log |Q(z)|\leq u_\phi^+(z) + {1\over 2} A_2\gamma \phi^{-1}\delta  \qquad \text{for}\quad z\in K.$$
On the other hand, we obtain as in the proof of Corollary \ref{cor_zeros_Laguerre} that the density of $\mu_\phi^+$ is bounded by a constant 
times $(1-\phi)^{-1}$. Therefore, since $\gamma$ is small, we deduce from
 the last assertion of Proposition \ref{prop_dist_sup} that 
 $$|\mu_{p,n}^+(I)-\mu_\phi^+(I)|\leq {\sqrt{\phi^{-1}\delta}\over 1-\phi}$$
for any interval $I$ in $\R$. Thus,
 $$|\mu_{p,n}(I)-\mu_\phi(I)|\leq {\sqrt{\phi\delta}\over 1-\phi}\leq {\sqrt{\delta}\over 1-\phi}$$
which ends the proof of the theorem (we just need to replace the constants $A_1,A_2$ in the theorem by $A_1',A_2'$).
\hfill $\square$


\small

\begin{thebibliography}{99}


\bibitem{BHZ}
Bai Z. D., Hu J., Zhou W., 
Convergence rates to the Marchenko-Pastur type distribution. 
{\it Stochastic Process. Appl.} {\bf122} (2012), no. 1, 68-92. 

\bibitem{BS}
Bai Z. D.,  Silverstein J. W.,
{\it Spectral analysis of large dimensional random matrices. }
Second edition. Springer Series in Statistics. Springer, New York, 2010.


\bibitem{BY}
Bai Z. D., Yin Y. Q., Convergence to the semicircle law. {\it Ann. Probab.} {\bf 16} (1988), no. 2, 863-875. 



\bibitem{BP}
Ben Arous G., P\'ech\'e S., 
Universality of local eigenvalue statistics for some sample covariance matrices.
{\it Comm. Pure Appl. Math.} {\bf 58} (2005), no. 10, 1316-1357. 

\bibitem{BEK}
Bloemendal A., Erd\"os L., Knowles A., Yau H.-T., Yin J., Isotropic local laws for sample covariance and generalized Wigner matrices. {\it Electron. J. Probab.} {\bf 19} (2014), no. 33, 53 pp.


\bibitem{CMS}
Cacciapuoti C., Maltsev A., Schlein B.,
Local Marchenko-Pastur law at the hard edge of sample covariance matrices. 
{\it J. Math. Phys.} {\bf 54} (2013), no. 4, 043302, 13 pp. 


\bibitem{Comtet}
Comtet L.,
{\it Advanced combinatorics. }
The art of finite and infinite expansions. Revised and enlarged edition. D. Reidel Publishing Co., Dordrecht, 1974.

\bibitem{DIW}
Dai D., Ismail M. E. H., Wang J.,
Asymptotics for Laguerre polynomials with large order and parameters.
{\it J. Approx. Theory}  {\bf 193} (2015), 4-19. 

\bibitem{DetteStudden}
 Dette H., Studden W. J., Some new asymptotic properties for the zeros of Jacobi, Laguerre, and Hermite polynomials. {\it Constr. Approx.} {\bf 11} (1995), no. 2, 227-238.


\bibitem{Dinh}
Dinh T.-C., Large deviation theorem for zeros of polynomials and Hermitian random  matrices. {\it Preprint}, 2016. 
{\tt  arXiv:1611.04271 }

\bibitem{DinhSibony}
 Dinh T.-C., Sibony N.,  Dynamics in several complex variables: 
 endomorphisms of projective spaces and polynomial-like mappings. 
 {\it Holomorphic dynamical systems,} 165-294, 
 Lecture Notes in Math., {\bf  1998}. {\it Springer, Berlin,} 2010.

\bibitem{GT}
 G\"otze F., Tikhomirov A., Rate of convergence in probability to the Marchenko-Pastur law. {\it Bernoulli} {\bf 10} (2004), no. 3, 503-548.

\bibitem{GS}
Grenander U., Silverstein J. W., Spectral analysis of networks with random topologies. {\it SIAM J. Appl. Math.} {\bf 32(2)} (1977), 499-519.

\bibitem{Jonsson}
Jonsson D.,  Some limit theorems for the eigenvalues of a sample covariance
matrix. {\it J. Multivariate Anal.} {\bf 12} (1982), 1-38.

\bibitem{MP}
Mar\v{c}enko V. A., Pastur L. A.,  Distribution for some sets of random
matrices. {\it Math. USSR-Sb.} {\bf 1} (1967), 457-483.

\bibitem{Szego}
 Szeg\"o G., Orthogonal polynomials. Fourth edition. {\it American Mathematical Society, Colloquium Publications,} Vol. {\bf XXIII}. 
 American Mathematical Society, Providence, R.I., 1975.


\bibitem{TaoVu}
Tao T., Vu V.,
Random covariance matrices: universality of local statistics of eigenvalues. {\it Ann. Probab.} {\bf 40} (2012), no. 3, 1285-1315.


\bibitem{Wachter}
Wachter K. W.,  The strong limits of random matrix spectra for sample matrices of independent elements. {\it Ann Probab.} {\bf 6(1)} (1978), 1-18.

\bibitem{Yin}
Yin Y. Q., Limiting spectral distribution for a class of random matrices. {\it J. Multivariate Anal.} {\bf 20} (1986), 50-68.

\end{thebibliography}
\end{document}